\documentclass[11pt]{article}
\usepackage[margin=1in]{geometry}
\usepackage{amsmath,amssymb,amsthm,booktabs}
\usepackage[hidelinks]{hyperref}
\hypersetup{pdftitle={Bundles with additional homotopy conditions},pdfauthor={A. V. Ershov}}
\newtheorem{theorem}{Theorem}[section]
\newtheorem{proposition}[theorem]{Proposition}
\newtheorem{lemma}[theorem]{Lemma}
\newtheorem{corollary}[theorem]{Corollary}
\theoremstyle{definition}
\newtheorem{definition}[theorem]{Definition}

\newcommand{\SU}{\mathrm{SU}}
\newcommand{\PU}{\mathrm{PU}}

\newcommand{\ZZ}{\mathbb Z}

\newcommand{\CC}{\mathbb C}

\newcommand{\RP}{\mathbb{RP}}

\newcommand{\Fr}{\mathrm{Fr}}
\newcommand{\Map}{\operatorname{Map}}
\newcommand{\Ad}{\operatorname{Ad}}
\newcommand{\End}{\operatorname{End}}
\newcommand{\holim}{\operatorname*{holim}}

\newcommand{\su}{\mathfrak{su}}
\newcommand{\Dalg}{\mathcal A}
\newcommand{\labspace}{\mathcal B}

\newcommand{\trivbundle}[1]{\underline{#1}}
\title{Bundles with additional homotopy conditions\\[3pt]
\large Biquotients, matrix Grassmannians and simplex dimension-drop algebras}
\author{A.\,V. Ershov}
\date{Revised manuscript}
\begin{document}
\maketitle
\begin{abstract}
A homotopy between two amplified special unitary bundles and a trivialization of a tensor product give two different finite-dimensional classifying spaces. For coprime ranks $k,l>1$ the first space is a biquotient of $\SU(kl)$ and the second is a homogeneous matrix Grassmannian. We distinguish the two manifolds for every such pair by their first Pontryagin classes. For $(k,l)=(2,3)$ their eighth integral cohomology groups have orders $3$ and $21$, respectively. The same amplification maps define a dimension-drop algebra on a labelled simplex. Its skeletal filtration gives a weighted cellular complex and computes $K_0$, $K_1$ and the unit when the vertex labels have greatest common divisor one. We also compute both groups for triangles with arbitrary vertex labels.
\end{abstract}

\section{Introduction}
Let $V$ and $W$ have complex ranks $k$ and $l$. A homotopy between the classifying maps of $V^{\oplus l}$ and $W^{\oplus k}$ gives one way of relating these bundles after increasing their ranks. A trivialization of $V\otimes W$ gives another. Both constructions use the same two commuting subgroups of $\SU(kl)$. In the first construction they act on opposite sides of $\SU(kl)$, and in the second they act on the same side. This distinction affects the topology of the classifying manifolds even when their loop spaces are homeomorphic as spaces.

The first construction arose in the theory of simplicial bundles in \cite{old,journal}. The second occurs in the study of algebra bundles embedded in a trivial matrix algebra bundle \cite{matrix,obstructions}. For coprime $k,l>1$ write $D_{k,l}$ for the double quotient and $M_{k,l}$ for the homogeneous matrix Grassmannian. Remark~32 of \cite{old} asked whether they could be homeomorphic. We obtain the following answers.
\begin{enumerate}
\item Both spaces have $H^4(-;\ZZ)\cong\ZZ$. The absolute coefficients of the first Pontryagin class in primitive generators are
\[
 2(k^2l^2-k^2-l^2) \quad\text{and}\quad 2|k^2-l^2|
\]
for $D_{k,l}$ and $M_{k,l}$, respectively. They distinguish the manifolds up to homeomorphism (Theorem~\ref{thm:pone}).
\item For $(k,l)=(2,3)$ the distinction persists up to homotopy. The degree-eight integral cohomology groups are $\ZZ/3$ for the biquotient and $\ZZ/21$ for the homogeneous space (Theorem~\ref{thm:eight}).
\item For the simplex algebra $\Dalg(d_0,\ldots,d_n)$, the first differential of its skeletal $K$-theory spectral sequence is exterior multiplication by the vector of vertex labels. If $\gcd(d_0,\ldots,d_n)=1$, the unital inclusion $\CC\to\Dalg(d_0,\ldots,d_n)$ induces an isomorphism on $K$-theory (Theorem~\ref{thm:k}). For a triangle we give the groups for every value of the gcd.
\end{enumerate}
The finite-stage bundle constructions also give an embedding criterion over complexes of dimension at most four. Section~\ref{sec:embed} makes explicit how it depends on the choice of a special unitary lift. Section~\ref{sec:comparison} describes the scope of the realization, tensor-product and extension constructions proposed in \cite{old}.

All classifying spaces use functorial models. Homotopy pullbacks and homotopy fibers are taken in spaces. Bases for bundle classification have CW type and the bundles are numerable. Products and mapping spaces can be taken in compactly generated spaces. We write $M_r$ for $M_r(\CC)$ and $I_r$ for its identity. Cohomology has integral coefficients unless a different coefficient group is indicated.

\section{Families on a labelled simplex}\label{sec:simp}
Let $[n]=\{0,\ldots,n\}$, with its usual order, and give its vertices labels $d_i>1$. For a nonempty subset $S$ put
\[
 E_S=\bigotimes_{i\in S}\CC^{d_i},\qquad d_S=\dim E_S,\qquad G_S=\SU(E_S).
\]
For $S\subset T$ the map $i_{ST}:G_S\to G_T$ inserts identity operators in the missing slots. The ordering and a fixed tensor-coordinate convention make these maps compatible under composition. We use the same notation for the analogous homomorphisms of matrix algebras. Unitary groups can be used throughout this section to describe complex bundles with unitary structure. Write $\Delta^S$ for the face on $S$, of dimension $|S|-1$.

\begin{definition}
A simplicial bundle of type $(d_0,\ldots,d_n)$ over $X$ is represented by maps
\[
 F_S:X\times\Delta^S\longrightarrow BG_S\qquad(\emptyset\ne S\subseteq[n])
\]
with the following compatibility. On $X\times\Delta^S$ the restriction of $F_T$ is $Bi_{ST}\circ F_S$ whenever $S\subset T$. Two families are concordant if there is such a family on $X\times[0,1]$ with the prescribed end families.
\end{definition}
Thus the vertices describe bundles. An edge gives a homotopy between their amplified classifying maps. On a triangle the new datum fills the loop made from the three edge homotopies after all three have been amplified into $BG_{012}$, where a string of vertex indices denotes their set. In bundle language concordances allow the endpoint bundles to vary through isomorphisms. Holding chosen endpoint representatives fixed gives a relative problem.

Let $\labspace(d_0,\ldots,d_n)$ denote the space of compatible families on a point, topologized as a subspace of the finite product of the face mapping spaces. The exponential law identifies concordance classes over $X$ with $[X,\labspace(d_0,\ldots,d_n)]$. There is also a homotopy-limit description
\begin{equation}\label{eq:faceholim}
 \labspace(d_0,\ldots,d_n)\simeq
    \holim_{\emptyset\ne S\subseteq[n]}BG_S
\end{equation}
with arrows going from smaller faces to larger faces. Indeed the standard homotopy-limit model consists of compatible maps from the nerve of the category of faces over $S$ into $BG_S$. This nerve is the barycentric subdivision of $\Delta^S$, naturally with respect to face inclusions \cite{BK}. This description uses no assumption about a choice of embeddings of the classifying spaces into a common ambient space.

For example, for three labels let $C_\partial$ denote the analogous space using only vertices and edges, and put $G=G_{012}$. Amplification of the boundary data defines a map $C_\partial\to\Map(S^1,BG)$. Then the full space is the following pullback
\begin{equation}\label{eq:triangle}
 C_\partial\times_{\Map(S^1,BG)}\Map(D^2,BG).
\end{equation}
Restriction from the disk is a fibration since the boundary inclusion is a cofibration. Therefore this is also a homotopy pullback. Every boundary loop has a filling since $BG$ is simply connected. Over the constant boundary the fiber is $\Omega^2BG\simeq\Omega G$. The boundary space itself fibers over $\prod_i BG_i$ with fiber of homotopy type $\prod_{i<j}G_{ij}$ over the basepoint tuple, by simultaneous endpoint evaluation of the three edge paths. Formula~\eqref{eq:triangle} records the filling and all endpoint compatibility in one diagram.

There is a finite model for the boundary when the three labels are pairwise coprime. Orient each edge from $i$ to $j$ for $i<j$. Put $H_\partial=\prod_i G_i$ and $E_\partial=\prod_{i<j}G_{ij}$. On the $ij$ entry the action of $(a_i)_i$ is
\[
 g_{ij}\longmapsto i_{i,ij}(a_i)\,g_{ij}\,i_{j,ij}(a_j)^{-1}.
\]
Here $i_{i,ij}$ abbreviates the inclusion for the corresponding singleton and edge. The boundary data are modeled by $EH_\partial\times_{H_\partial}E_\partial$: this is the bundle over $\prod_i BG_i$ of edge identifications between the two amplified vertex fibers for each edge. It is obtained by taking the three associated bundles of identifications over that common base and forming their fiber product there. The freeness test in the next theorem, applied on each edge, forces each vertex stabilizer coordinate to be the identity. Thus $C_\partial\simeq E_\partial/H_\partial$ in this case. Formula~\eqref{eq:triangle} then describes the remaining filling problem over a finite-dimensional boundary model.

\section{Two classifying manifolds}\label{sec:two}
Fix $k,l>1$. Put $G=\SU(kl)$, $H=\SU(k)\times\SU(l)$, and write
\[
 \alpha(a)=a\otimes I_l,\qquad \beta(b)=I_k\otimes b,\qquad
 \tau(a,b)=a\otimes b .
\]
A one-simplex with labels $k,l$ gives the homotopy pullback
\begin{equation}\label{eq:double}
 B\SU(k)\times^{h}_{BG}B\SU(l)
\end{equation}
of $B\alpha$ and $B\beta$. A second space is the homotopy fiber of $B\tau:BH\to BG$ over the basepoint. In the first space the amplified bundles are identified. In the second their tensor product has a specified special unitary trivialization. In either case families are taken up to concordance, including the identifications on the ends.

\begin{theorem}[Finite quotient models {\cite{old}}]\label{thm:models}
If $\gcd(k,l)=1$, the two spaces just defined are modeled respectively by the smooth compact manifolds
\[
 D_{k,l}=\alpha(\SU(k))\backslash G/\beta(\SU(l)),\qquad
 M_{k,l}=G/\tau(H).
\]
The second manifold parametrizes unital $*$-subalgebras of $M_{kl}$ isomorphic to $M_k$. The associated rank-$k$ and rank-$l$ bundles $V,W$ on each quotient satisfy, respectively,
\begin{equation}\label{eq:twoid}
 V^{\oplus l}\cong W^{\oplus k}\quad(D_{k,l}),\qquad
 V\otimes W\cong\trivbundle{\CC}^{\,kl}\quad(M_{k,l}).
\end{equation}
\end{theorem}
\begin{proof}
For the action $(a,b)g=\alpha(a)g\beta(b)^{-1}$, the Borel construction $EH\times_HG$ is a model of \eqref{eq:double}. One way to see this is to pull back the diagonal model $G\to B\Delta G\to BG\times BG$ along $(B\alpha,B\beta):BH\to BG\times BG$. The fiber $G$ parametrizes identifications of the two extended fibers. If $\alpha(a)=g\beta(b)g^{-1}$, each eigenvalue has multiplicity divisible by both $l$ and $k$. Coprimality and the dimension $kl$ force a scalar. The determinant conditions put this scalar in $\mu_k\cap\mu_l=\{1\}$, where $\mu_r$ denotes the $r$th roots of unity. The action is free. The projection of the Borel construction onto the orbit space has contractible fiber $EH$ and is an equivalence. The tautological identification on this quotient is given by $g$ itself, which intertwines the right and left associated fibers.

The tensor homomorphism $\tau$ has kernel given by the opposite scalars in $\mu_k\cap\mu_l$. It therefore embeds $H$ as a closed subgroup in the coprime case. The standard fiber sequence $G/\tau(H)\to BH\to BG$ proves the second assertion. The tensor representation extends to the standard representation of $G$, giving the indicated trivialization. Finally, conjugation by $\PU(kl)$ is transitive on the unital copies of $M_k$, with stabilizer $\PU(k)\times\PU(l)$ in tensor position. This gives the same quotient since $\mu_k\mu_l=\mu_{kl}$. The second matrix algebra is the centralizer of the first.
\end{proof}

\paragraph{Path groups.}
For closed subgroups $K,L\subset G$ let $P_{K,L}G$ be the group of paths in $G$ with first value in $K$ and last value in $L$, under pointwise multiplication. Write $K=\alpha(\SU(k))$, $L=\beta(\SU(l))$ and $J=\tau(H)$. The groups $P_{K,L}G$ and $P_{J,1}G$ have $D_{k,l}$ and $M_{k,l}$, respectively, as classifying spaces up to homotopy. To see this, pull the endpoint homomorphism $G^{[0,1]}\to G^2$ back to $K\times L$, or to $J\times\{1\}$. It is onto with local sections since $G$ is a connected Lie group, and its kernel is $\Omega G$. The original and pulled-back group extensions give fiber sequences on classifying spaces with the same fiber $B\Omega G$. Their square is therefore homotopy Cartesian. The inclusion of constant paths $G\to G^{[0,1]}$ is a homotopy equivalence of groups and induces the diagonal on endpoint values. This identifies the two classifying spaces with the homotopy pullback and the homotopy fiber already computed. The spaces of paths are homeomorphic by
\[
 \gamma\longmapsto\big(t\longmapsto\gamma(t)\gamma(1)^{-1}\big).
\]
Indeed, for the inverse factor the initial value of $\kappa\in P_{J,1}G$ uniquely as $a b^{-1}$ with $a\in K$, $b\in L$, and set $\gamma(t)=\kappa(t)b$. Factorization is continuous since $K\times L\to J$, $(a,b)\mapsto ab^{-1}$, is a homeomorphism. The formula changes the pointwise multiplication by conjugations involving the endpoint. The next section detects the effect on the classifying spaces.

\section{Distinguishing the manifolds}\label{sec:dist}
We continue to assume $k,l>1$ and $\gcd(k,l)=1$. On either manifold set $x=c_2(V)$ and $y=c_2(W)$.

\begin{lemma}\label{lem:degreefour}
Both manifolds are $3$-connected and have $H^4(-;\ZZ)\cong\ZZ$. With a suitable primitive generator $z$ on each space,
\[
 (x,y)=(kz,lz)\quad\text{on }D_{k,l},\qquad
 (x,y)=(kz,-lz)\quad\text{on }M_{k,l}.
\]
\end{lemma}
\begin{proof}
In either principal $H$-bundle with total space $G$, the map on $\pi_3$ induced by the orbit through the identity has coefficients $(l,-k)$ for the two-sided action and $(l,k)$ for tensor multiplication, up to an overall sign. Indeed a block sum of $s$ copies acts on $\pi_3\SU$ by multiplication by $s$, and inversion acts by negation. These statements follow equally by clutching over $S^4$ and applying $c_2$. The maps are onto. As $H,G$ are simply connected with vanishing $\pi_2$ and $\pi_4(G)=0$ by Bott periodicity, the exact sequence gives $3$-connectivity and identifies $\pi_4$ of the base with the integral kernel of this map. The kernel is generated by $(k,l)$ or $(k,-l)$, respectively. The two entries are the clutching numbers of the two bundles pulled back to a generating four-sphere. Hurewicz and the universal coefficient theorem now give the claimed primitive normalization.
\end{proof}

\begin{theorem}\label{thm:pone}
For every coprime pair $k,l>1$ the manifolds $D_{k,l}$ and $M_{k,l}$ are not homeomorphic.
\end{theorem}
\begin{proof}
Consider first a general free action of a compact group $H$ on $G$ by $h\cdot g=\lambda(h)g\rho(h)^{-1}$ and let $Q$ be the quotient. The vertical tangent bundle of $G\to Q$ is equivariantly $G\times\mathfrak h$ with the adjoint action on $\mathfrak h$. Left trivialization of $TG$ identifies the action on the Lie-algebra coordinate with $\Ad\rho(h)$. Splitting the vertical tangent sequence by an invariant metric and descending to $Q$ yields
\begin{equation}\label{eq:tangent}
 TQ\oplus(G\times_H\mathfrak h)\cong G\times_H\mathfrak g_{\Ad\rho}.
\end{equation}
This is the stable tangent formula for biquotients, cf.\ \cite{Singhof}. On $D_{k,l}$ the bundle on the right is $\su(U)$, where $U=W^{\oplus k}\cong V^{\oplus l}$. On the homogeneous quotient it is trivial. For instance, use the right action by $\tau(H)$ and the map $[g,X]\mapsto(g\tau(H),\Ad_g X)$ in the associated-bundle convention.

For a rank-$r$ special unitary bundle $E$, the complexification of $\su(E)$ is $\End_0(E)$. If $s_i$ are the Chern roots of $E$, the degree-four component of the Chern character of $\End_0(E)$ is $\tfrac12\sum_{i,j}(s_i-s_j)^2=-2r c_2(E)$. Equivalently the splitting principle gives the integral identity $p_1(\su(E))=-2r c_2(E)$. The vertical term in \eqref{eq:tangent} is $\su(V)\oplus\su(W)$. Since $c_2(U)=lx=ky$ on $D_{k,l}$, we get
\begin{align*}
 p_1(TD_{k,l})&=-2kl(lx)+2kx+2ly,\\
 p_1(TM_{k,l})&=2kx+2ly .
\end{align*}
By Lemma~\ref{lem:degreefour} these classes have absolute coefficients $2(k^2l^2-k^2-l^2)$ and $2|k^2-l^2|$ in integral generators. The former exceeds the latter since their difference is $2\max(k^2,l^2)(\min(k^2,l^2)-2)>0$. A homeomorphism would act by a sign on the integral groups and preserve the rational Pontryagin class by Novikov's theorem \cite{Novikov}. This is impossible.
\end{proof}

\begin{theorem}\label{thm:eight}
The first coprime pair already gives distinct homotopy types. More precisely,
\[
 H^8(D_{2,3};\ZZ)=\ZZ/3,\qquad H^8(M_{2,3};\ZZ)=\ZZ/21.
\]
\end{theorem}
\begin{proof}
Use the fibration $G\to EH\times_H G\to BH$ with $G=\SU(6)$ and $H=\SU(2)\times\SU(3)$ for each of the two free actions. Write $a=c_2(V_0)$, $b=c_2(W_0)$ and $c=c_3(W_0)$ for the universal factors on $BH$. The base ring is $\ZZ[a,b,c]$, with $a,b$ of degree four and $c$ of degree six. The fiber has exterior generators $u_3,u_5,u_7,u_9,u_{11}$. In the integral Serre spectral sequence the transgression of $u_{2i-1}$ is the difference of the $i$th Chern classes of the two representations, taken modulo earlier differentials. This follows by pulling back the universal diagonal fibration along the two representations \cite[Theorem~1]{Eschenburg}. For the homogeneous action one representation is trivial, so we use the Chern classes of $V_0\otimes W_0$, choosing signs for the generators.

In degree four the relations are $3a-2b$ for $D$ and $3a+2b$ for $M$. The next relation affecting the base in total degree eight is obtained from the fourth Chern classes
\begin{align*}
 c_4(V_0^{\oplus3})-c_4(W_0^{\oplus2})&=3a^2-b^2,\\
 c_4(V_0\otimes W_0)&=3a^2+b^2 .
\end{align*}
For the second equality let $x,-x$ be the roots of $V_0$, so that $a=-x^2$. The roots of $W_0$ satisfy $\sum y_i=0$, $\sum_{i<j}y_i y_j=b$, $y_1y_2y_3=c$. The Chern polynomial with leading coefficient one is
\[
 \prod_{i=1}^3(t+y_i+x)(t+y_i-x)
  =\big((t+x)^3+b(t+x)+c\big)
   \big((t-x)^3+b(t-x)+c\big).
\]
The coefficient of $t^2$ is $3x^4+b^2$. The first equality follows by multiplying the total Chern classes of the direct sums.

In bidegree $(8,0)$ the $d_4$ images are the multiples of $3a\mp2b$ by $a$ and $b$. The $d_6$ into this bidegree has zero source since $H^2(BH)=0$. Modulo the degree-four relation the degree-eight part is freely generated by $z^2$, with $a=2z$ and $b=3z$ for $D$ or $b=-3z$ for $M$. Thus $d_8(u_7)$ imposes $3z^2=0$ or $21z^2=0$. There are no later differentials into $(8,0)$. The only other nonzero $E_2$ term of total degree eight is generated by $u_3u_5$ in $(0,8)$. Its $d_4$ is $(3a\mp2b)u_5\ne0$, so it contributes no further filtration quotient. This proves the exact orders as asserted.
\end{proof}

The same pair has first Pontryagin coefficients of absolute values $46$ and $10$. Theorem~\ref{thm:pone} answers the homeomorphism question for all coprime pairs. The stronger distinction by integral cohomology in Theorem~\ref{thm:eight} concerns the indicated pair.

\section{Lifts and embeddings}\label{sec:embed}
\paragraph{Lifting a prescribed bundle.}
For $2\le r<\min(k,l)$ consider a map $f:S^{2r}\to B\SU(kl)$. Through the homotopy pullback model there is a comparison $D_{k,l}\to B\SU(kl)$, unique up to the specified homotopy between the two projections followed by amplification. Bott periodicity identifies the groups in degree $2r$ of the three special unitary classifying spaces with $\ZZ$. Amplification from the two factors acts by $l$ and $k$, respectively. Thus $f$ lifts up to homotopy precisely when
\[
  [f]\in l\ZZ\cap k\ZZ=kl\ZZ.
\]
Indeed belonging to the intersection supplies maps to the two factors and homotopies from each amplification to $f$, hence supplies the path in the pullback. Extending an equatorial map while allowing the ambient bundle to vary is a different lifting condition. For $k=3,l=4,r=2$ every map $S^3\to D_{3,4}$ extends over a hemisphere by Lemma~\ref{lem:degreefour}, while a generator on $S^4$ fails the displayed test.

\paragraph{Embedding a matrix algebra bundle.}
Put $\Fr_{k,l}=\operatorname{Emb}_1^*(M_k,M_{kl})$. Conjugacy of unital representations gives
\[
 \Fr_{k,l}=\PU(kl)/(1_k\otimes\PU(l)).
\]
The group $\PU(k)$ acts freely by changes of source coordinates and the orbit space is $M_{k,l}$ for coprime ranks. For a principal $\PU(k)$-bundle $P$ on $X$ an embedding of the associated $M_k$-bundle in $X\times M_{kl}$ is a section of $P\times_{\PU(k)}\Fr_{k,l}$. Applying this description to the universal bundle yields the fiber sequence
\begin{equation}\label{eq:embed}
 \Fr_{k,l}\longrightarrow M_{k,l}\longrightarrow B\PU(k),
\end{equation}
where the middle term is identified with the Borel construction by freeness. The universal reduction through $B\SU(k)$ follows from the rank-$k$ factor on $M_{k,l}$. Hence an embedding forces the obstruction to lifting $P$ through $\SU(k)\to\PU(k)$ to vanish in $H^2(X;\ZZ/k)$.

\begin{proposition}\label{prop:obstruction}
Let $\dim X\le4$, where $X$ is a finite CW complex, and fix coprime $k,l>1$. A projective rank-$k$ bundle admits a fiberwise embedding of its matrix algebra into the trivial $M_{kl}$-bundle if and only if it admits a special unitary lift $V$ with $c_2(V)\in kH^4(X;\ZZ)$.
\end{proposition}
\begin{proof}
An embedding pulls back the two universal factors, with $V\otimes W$ trivial. Since $c_1(V)=c_1(W)=0$, the identity $lc_2(V)+kc_2(W)=0$ holds. Reducing in $H^4(X;\ZZ)/kH^4(X;\ZZ)$ and using coprimality proves necessity. Conversely if $c_2(V)=ku$, choose a rank-$l$ special unitary bundle $W$ with $c_2(W)=-lu$. Indeed $c_2:B\SU(s)\to K(\ZZ,4)$ induces an isomorphism through degree four for each $s\ge2$, so $c_2$ classifies such bundles on $X$. The tensor product has vanishing $c_2$, hence admits a special unitary trivialization. This trivialization embeds $\End(V)$ in $X\times M_{kl}$.
\end{proof}

For a fixed lift the condition is equivalent to $c_2(V)\bmod k=0$ by the coefficient sequence. Changing the lift may change this class. On $\RP^4$ let $L$ be the complex line of the sign representation. Its square is trivial, and both $\trivbundle{\CC}^{\,2}$ and $L\oplus L$ lift the same trivial projective bundle. If $t$ generates $H^1(\RP^4;\ZZ/2)$, the mod-two reduction of $c_1(L)$ is $t^2$, since the underlying real plane is two copies of the real sign line. Thus $c_2(L\oplus L)$ reduces to $t^4\ne0$, whereas the trivial lift has zero second Chern class. The existential quantifier in Proposition~\ref{prop:obstruction} accounts for this indeterminacy. Higher-dimensional lifting problems retain further obstructions \cite{obstructions}.

\section{Dimension-drop algebras on a simplex}\label{sec:k}
The identity insertions of Section~\ref{sec:simp} can be applied to functions with values in matrix algebras. We use the minimal tensor product throughout. The relation of these boundary conditions to operator algebras was pointed out by Ralph Meyer in connection with the original work \cite{old}.

\begin{definition}
Write $D=\prod_{i=0}^n d_i$. The simplex dimension-drop algebra $\Dalg(d_0,\ldots,d_n)$ is the unital $C^*$-subalgebra of $C(\Delta^n,M_D)$ consisting of functions whose restriction to each nonempty face $\Delta^S$ takes values in the standard copy of $\End(E_S)$, amplified by identities on the complementary tensor factors.
\end{definition}
For one edge this is the dimension-drop interval algebra
\[
 \Dalg(k,l)=\{f\in C([0,1],M_k\otimes M_l):
 f(0)\in M_k\otimes I_l,\quad f(1)\in I_k\otimes M_l\}.
\]
For $\gcd(k,l)=1$ these are the prime dimension-drop building blocks used by Jiang and Su \cite{JS}. The edge path group $P_{K,L}G$ acts by pointwise conjugation. A principal bundle with this group therefore gives a bundle of these algebras together with its indicated structure group. The endpoint evaluations are equivariant for conjugation on the corresponding matrix factors.

\begin{proposition}[The two endpoint maps]\label{prop:interval}
Let $g=\gcd(k,l)$. In rank-one projection coordinates the image of $K_0(\Dalg(k,l))$ in $\ZZ^2$ is generated by $(k/g,l/g)$, its unit is $(k,l)$, and $K_1(\Dalg(k,l))=\ZZ/g$. For $g=1$, the ideal $J=\ker(\mathrm{ev}_0:\Dalg(k,l)\to M_k)$ has $K_0(J)=0$ and $K_1(J)=\ZZ/k$.
\end{proposition}
\begin{proof}
Simultaneous endpoint evaluation is onto, with kernel $C_0((0,1),M_{kl})$. The six-term sequence identifies $K_0(\Dalg(k,l))$ and $K_1(\Dalg(k,l))$ with the kernel and cokernel of
\[
 \ZZ^2\longrightarrow\ZZ,\qquad (a,b)\longmapsto la-kb .
\]
Indeed the extension is the pullback of the two-endpoint path extension of $M_{kl}$ under amplification at the two ends. Its connecting map is the difference of their $K_0$ maps, up to the choice of interval orientation. This proves the first assertions. For $g=1$ the unit generates $K_0(\Dalg(k,l))$ and evaluation at zero sends it to $k$ in rank-one coordinates. The six-term sequence for this single endpoint now gives the asserted groups of $J$.
\end{proof}

\begin{theorem}[Skeletal computation]\label{thm:k}
Let $g=\gcd(d_0,\ldots,d_n)$ and $E=\ZZ^{n+1}$ with basis $e_0,\ldots,e_n$. On each even coefficient row the first differential in the skeletal $K$-theory spectral sequence for $\Dalg(d_0,\ldots,d_n)$ is the differential of the complex
\begin{equation}\label{eq:koszul}
 \bigwedge^1E\xrightarrow{\,d\wedge\,}\bigwedge^2 E
 \xrightarrow{\,d\wedge\,}\cdots\xrightarrow{\,d\wedge\,}\bigwedge^{n+1}E,
 \qquad d=\sum_{i=0}^n d_i e_i .
\end{equation}
The leftmost term is in cellular degree zero. On these rows the second page is
\[
 E_2^{p,q}\cong
 \begin{cases}
 \ZZ(d/g),&p=0,\\
 (\ZZ/g)^{\binom np},&1\le p\le n.
 \end{cases}
\]
The odd coefficient rows vanish. In particular, if $g=1$ then
\[
 K_0(\Dalg(d_0,\ldots,d_n))=\ZZ[1],\qquad
 K_1(\Dalg(d_0,\ldots,d_n))=0.
\]
\end{theorem}
\begin{proof}
Let $A_p$ denote the algebra of compatible sections over the $p$-skeleton, with $A_{-1}=0$. Restriction $A_p\to A_{p-1}$ is onto because the matrix entries of the given section on the boundary of each $p$-face extend over that face in its own matrix algebra. These extensions agree on intersections because distinct open $p$-faces meet only along their boundaries after closure. The kernel is
\begin{equation}\label{eq:cells}
 \bigoplus_{|S|=p+1} C_0(\mathring{\Delta}^S)\otimes\End(E_S).
\end{equation}
Consequently the kernels of restriction from $A_n$ give a finite filtration by closed ideals. The exact couples of the six-term sequences give the spectral sequence of this filtration \cite{Schochet}. In cohomological indexing the $S$-summand of $E_1^{p,q}$ is $K_{p+q}(C_0(\mathring{\Delta}^S)\otimes\End(E_S))=K_q(\End(E_S))$ by suspension and periodicity. The abutment in total degree $p+q$ is $K_{p+q}(A_n)$, with the $K$-degrees taken modulo two.

Choose the ordered vertex orientations on all faces. To compute the incidence entry from $S$ to $T$ with $|T|=|S|+1$, restrict the connecting map to the closed cell $\Delta^T$. If $S$ is a facet, extending by zero from $\mathring\Delta^S$ to $\partial\Delta^T$ and taking the disk boundary map gives the usual signed incidence number of that facet, by the suspension identification. On coefficients the restriction uses $i_{ST}:\End(E_S)\to\End(E_T)$. A rank-one projection is taken to a projection of rank $d_j$, where $T=S\cup\{j\}$. The combined entry is therefore $(-1)^{|\{i\in S:i<j\}|}d_j$ with a consistent choice of the suspension signs. All other entries are zero. Identifying the $S$-generator with the ordered wedge of its $e_i$ gives \eqref{eq:koszul}.

The primitive vector $v=d/g$ extends to a basis of $E$; write $E=\ZZ v\oplus F$, $\operatorname{rank}F=n$. In degree $p$ the summand $\bigwedge^{p+1}F$ maps injectively to $v\wedge\bigwedge^{p+1}F$ in the next degree by multiplication by $g$ (with a zero source at the top). The summand $v\wedge\bigwedge^p F$ maps to zero. In degree zero its generator $v$ has no incoming differential. In positive degrees the boundaries in this summand are precisely its multiples by $g$. This proves the page calculation.

For $g=1$ the page is concentrated in the zero column, so the finite filtration has neither further differentials nor additive extensions. Under the edge map to the vertices the unit has rank vector $(d_0,\ldots,d_n)=d$. This is the generator just found. This also identifies the unital map from the scalars on both $K$-groups.
\end{proof}

\begin{corollary}[A triangle with arbitrary labels]\label{cor:triangle}
If $g=\gcd(d_0,d_1,d_2)$, then as abelian groups
\[
 K_0(\Dalg(d_0,d_1,d_2))\cong\ZZ\oplus\ZZ/g,\qquad
 K_1(\Dalg(d_0,d_1,d_2))\cong(\ZZ/g)^2 .
\]
The image of the unit in the free quotient is $g$ times its primitive positive-rank generator.
\end{corollary}
\begin{proof}
The cellular degrees are $0,1,2$. The only possible higher differential by column length is $d_2$, which changes coefficient parity and vanishes. Column one gives all of $K_1$. In even total degree the filtration gives an extension of the free group in column zero by the cyclic group in column two. It splits as a group since its quotient is free. The evaluation at the vertices sends the unit to $d=g v$, proving the assertion about its image in the quotient.
\end{proof}
For example, the labels $(6,10,15)$ give $K_0=\ZZ[1]$ and $K_1=0$, although each pair has a common factor. For $(2,4,6)$ one obtains $K_0\cong\ZZ\oplus\ZZ/2$ and $K_1\cong(\ZZ/2)^2$. For higher simplices with $g>1$, Theorem~\ref{thm:k} supplies a second page. The possible odd-length differentials and the additive extension problems must be addressed before reading off integral groups from that page.

\section{Stabilization and related constructions}\label{sec:comparison}
\subsection*{A common stable label}
There is a precise sense in which the classifying spaces of Section~\ref{sec:two} approximate the same space. Write $B=\mathrm{BSU}_{\otimes}$ for $\mathrm{BSU}$ with multiplication on stable virtual bundles of rank one induced by tensor product. For a finite CW complex $X$, $[X,B]$ consists of the rank-one elements $\xi\in K^0(X)$ with $c_1(\xi)=0$. The rank-zero ideal is nilpotent, so these elements have multiplicative inverses.

\begin{proposition}\label{prop:stable}
There are comparisons $D_{k,l}\to B$ and $M_{k,l}\to B$, defined up to homotopy. Each induces isomorphisms on $\pi_i$ for $i\le 2\min(k,l)$. The universal rank-one labels have the following characterizations in $K^0$:
\[
\begin{array}{c|cc}
 & [V] & [W]\\ \hline
 D_{k,l}&k\xi&l\xi\\
 M_{k,l}&k\xi&l\xi^{-1}.
\end{array}
\]
\end{proposition}
\begin{proof}
On $D_{k,l}$ we have $l[V]=k[W]$. Choose $s,t$ with $sk+tl=1$. Then $\xi=s[V]+t[W]$ has rank one, $c_1=0$, and satisfies the two identities, which also determine it uniquely. On $M_{k,l}$ one has $[V][W]=kl$. In $K^0(M_{k,l})[1/k]$ take $[V]/k$, and in $K^0(M_{k,l})[1/l]$ take $([W]/l)^{-1}$. The latter inverse exists by nilpotence of the reduced ideal. These elements agree after inverting $kl$. Since $k,l$ are coprime, the sequence
\[
 0\longrightarrow A\longrightarrow A[1/k]\oplus A[1/l]
 \xrightarrow{\ \mathrm{difference}\ } A[1/(kl)]
\]
is exact for any abelian group $A$. For instance, matching fractions $a/k^r,b/l^s$ can, by increasing the powers, be taken to satisfy $l^s a=k^r b$ in $A$. A B\'ezout identity for $k^r,l^s$ then gives a common preimage. Thus the local classes yield a unique rank-one $\xi$ satisfying the claimed identities. Since $k c_1(\xi)=l c_1(\xi)=0$ they also give $c_1(\xi)=0$. On both spaces this class defines the comparison map.

For $4\le 2j\le 2\min(k,l)$ the orbit sequence identifies $\pi_{2j}$ with the kernel of $(l,-k)$ or $(l,k)$ on the two copies of $\ZZ=\pi_{2j-1}\SU$ for the factors. Indeed these factor groups are stable and $\pi_{2j}\SU(kl)=0$. The first projection of a kernel generator is $k$ times a generator. Since $[V]=k\xi$, each comparison sends that generator to a Bott generator of $\pi_{2j}B$. The odd groups in the stated range vanish by the surjectivity of the same orbit maps in odd degree and vanishing of the preceding even groups. The degrees at most three follow from Lemma~\ref{lem:degreefour}.
\end{proof}
In particular on finite complexes of dimension at most $2\min(k,l)$ the two comparisons classify rank-one labels by passage to the corresponding Postnikov section. With products given by tensoring the two vertex factors, the identity between amplifications, or the trivialization of the tensor product, also tensors. When the new pair of ranks is coprime, uniqueness in Proposition~\ref{prop:stable} gives the external tensor product of the two labels on the product of bases. A multiplicative tower such as $(2^r,3^r)$ therefore has compatible comparisons, and its mapping telescope maps by a weak equivalence to $B$ in either construction. One chooses homotopies on successive cylinders using equality of the stable labels. The increasing stable range then proves the equivalence on all homotopy groups. For coordinatewise divisibility the totality of coprime pairs is not directed: $(2,3)$ and $(3,2)$ have no common admissible multiple. For a localization or a torsion-coefficient limit one must specify the indexing diagram and its effect on homotopy groups.

\begin{corollary}[Changing the complementary rank]
Over a finite CW base $X$, suppose a bundle of algebras with fiber $M_k$ embeds unitally in $X\times M_{kl}$ for some $l$ coprime to $k$. Then it embeds unitally in $X\times M_{km}$ for every sufficiently large $m$.
\end{corollary}
\begin{proof}
The subalgebra family gives a map $X\to M_{k,l}$. Pull back its factors and its label $\xi$ from Proposition~\ref{prop:stable}; in particular the given algebra bundle is $\End(V)$ and $[V]=k\xi$. For sufficiently large $m$, the stable special unitary class $m\xi^{-1}$ of rank $m$ has a rank-$m$ representative $W_m$ by unitary stability. The bundle $V\otimes W_m$ is stably trivial and has rank $km$ in the cancellation range, hence is trivial. Acting on its first factor gives the embedding. For instance, $2m>\dim X$ suffices for these stability arguments when $m\ge2$.
\end{proof}

\subsection*{Coherence and the operator-algebra extensions}
The face restrictions of Section~\ref{sec:simp} define a semisimplicial space after taking the disjoint union over the label sequences in each degree. Subfamilies closed under deleting vertices give semisimplicial subspaces; the condition of pairwise coprimality alone has this property. A condition imposing a positive lower bound on simplex dimension is not face-closed. Identifying the realization of any such space with $B\SU$ or $B\PU$ requires a comparison for that realization. Proposition~\ref{prop:stable} gives a comparison for the specified telescopes. On a numerable cover the same boundary maps give coherent local data by putting the face corresponding to $S$ on the intersection $U_S$ and requiring agreement after restriction and amplification. Classification of the resulting descent data requires a specified notion of refinement and equivalence.

For tensor powers of the interval algebra the parameters are independent. The center of $\Dalg(k,l)^{\otimes m}$ is $C([0,1]^m)$, whereas the center of $\Dalg(k^m,l^m)$ is $C([0,1])$. Indeed the tensor product embeds in the continuous matrix-valued functions on the cube, its fibers on the open cube are full matrix algebras, and it contains all continuous scalar functions. Centrality on that dense open cube forces scalar values everywhere. For $m>1$ the centers distinguish these algebras. Restriction to the diagonal, with a permutation of tensor coordinates, does give a unital $*$-homomorphism to $\Dalg(k^m,l^m)$. Constructing a multiplicative infinite-loop-space model from such maps requires symmetric coherence and a group-completion comparison. Existing models of units for tensor product are provided by Segal and by Dadarlat--Pennig \cite{Segal,DP}.

Finally the cone $C_0([0,1),M_{kl})$ evaluates at zero onto $M_{kl}$. The inverse image of the proper unital subalgebra $I_k\otimes M_l$ under that evaluation is a subalgebra. It cannot be an ideal, since evaluation is onto and its image is not an ideal of $M_{kl}$. The endpoint ideal $J$ of Proposition~\ref{prop:interval} gives a genuine extension with quotient $M_k$. It is essential since it contains the open-interval ideal, whose fibers are full on a dense open set. Left and right multiplication therefore embed $\Dalg(k,l)$ in the multiplier algebra $\mathcal M(J)$ and induce a Busby map $M_k\to\mathcal M(J)/J$. Delooping a sequence of classifying spaces or classifying Brauer-type invariants would require a further construction from this extension.

\paragraph{Acknowledgments.}
I thank E.\,V. Troitsky for his sustained interest and support, and L.\,A. Alania, V.\,M. Manuilov and A.\,S. Mishchenko for discussions. Much of the original research was carried out during a visit to G\"ottingen supported by the RFFI--DFG program. I am grateful to Thomas Schick for his hospitality and our discussions there, and to Ralph Meyer for drawing my attention to the operator-algebra interpretation. The original work received support from RFFI grant 07-01-00046 and RFFI--DFG grant 07-01-91555.

\paragraph{Use of AI tools.}
OpenAI's ChatGPT assisted in reviewing the earlier proofs, developing and checking the calculations, and preparing this revised manuscript.

\end{document}